\documentclass[11pt]{article}
\usepackage[utf8]{inputenc}
\usepackage{amsmath}
\usepackage{fullpage}
\usepackage{commath}
\usepackage{amsfonts}
\usepackage{amssymb}
\usepackage{amsthm}
\usepackage{tikz}
\usepackage{caption}
\usepackage{placeins}
\newtheorem{thm}{Theorem}[section]
\newtheorem{lem}{Lemma}[section]

\newtheorem{cor}{Corollary}[section]
\newtheorem{rmk}{Remark}[section]
\newtheorem{prob}{Problem}[section]
\DeclareMathOperator{\diam}{diam}
\DeclareMathOperator{\rank}{rank}

\DeclareMathOperator{\In}{In}
\newcommand{\E}{\mathcal{E}}

\title{Irreducibility, inertia and invertibility of the eccentricity matrices of unicyclic graphs}

\author{%
  Subhendu Chakraborty\textsuperscript{1,2,}\thanks{\texttt{subhendu.chakraborty.math@gmail.com}},~
  Iswar Mahato\textsuperscript{3,}\thanks{\texttt{iswarmahato02@gmail.com}, \texttt{iswar.m@srmap.edu.in}}\\[2ex]
  {\small \textsuperscript{1}Institute for Advancing Intelligence (IAI), TCG CREST, Kolkata-700091, India}\\
  {\small \textsuperscript{2}Academy of Scientific and Innovative Research (AcSIR), Ghaziabad-201002, India}\\
  {\small \textsuperscript{3}Department of Mathematics, SRM University-AP, Andhra Pradesh-522240, India}%
}

\date{}
\begin{document}
\maketitle
\begin{abstract}
The eccentricity matrix of a connected graph $G$, denoted by $\mathcal{E}(G)$, is obtained from the distance matrix of $G$ by keeping each entry that is the largest in its row or the largest in its column, and replacing all other entries by $0$. The eigenvalues of $\mathcal{E}(G)$ are called the $\mathcal{E}$-eigenvalues of $G$. In this article, we consider the eccentricity matrix of unicyclic graphs and study its various spectral properties. First, we give two sufficient conditions for the irreducibility of the eccentricity matrix of a connected graph, which generalize some previously known results. Then, we show that the eccentricity matrix of every unicyclic graph of girth $3$, $5$ or $7$ is irreducible and provide a unicyclic graph whose eccentricity matrix is reducible for every odd girth $g\ge 9$. On the other hand, we also provide a class of unicyclic graphs of odd girth, whose eccentricity matrices are always irreducible. Additionally, we find the inertia of eccentricity matrices of unicyclic graphs in which the trees attached to the cycle are nontrivial and have equal heights. Moreover, we characterize the graphs in this class whose E-eigenvalues are symmetric about the origin. Finally, we obtain necessary conditions for the invertibility of the eccentricity matrix of a unicyclic graph.
\end{abstract}




\section{Introduction}\label{sec1}
Throughout this paper, we consider finite, simple and connected graphs. Let $G$ be a graph with vertex set $V(G)$ and edge set $E(G)$. The number of vertices in $G$ is the \textit{order} of $G$. The \textit{adjacency matrix} of a graph $G$ on $n$ vertices, denoted by $A(G)$, is the $n \times n$ matrix whose rows and columns are indexed by the vertex set of $G$ and the  $(u,v)$-th entry $A(G)_{uv}=1$, if the vertices $u$ and $v$ are adjacent and $A(G)_{uv}=0$ otherwise.
The \textit{distance} between the vertices $u,v\in V(G)$, denoted by $d_G(u,v)$, is the length of a shortest path between them in $G$, and $d_G(u,u)=0$ for all $u \in V(G)$. The \textit{distance matrix} of a graph $G$ of order $n$ is the $n \times n$ matrix with its rows and columns indexed by the vertices of $G$, and the $(u,v)$-th entry is equal to $d_G(u,v)$. The \textit{eccentricity} $e_G(u)$ of a vertex $u\in V(G)$ is defined as $e_G(u)=\max\{d_G(u,v):v\in V(G)\}$.  A graph is \emph{self-centered} if all its vertices have the same eccentricity. The \textit{diameter} of $G$, denoted by $\diam(G)$, is the maximum distance between any two vertices in $G$. A \textit{diametrical path} in $G$ is a shortest path between two vertices $u$ and $v$ of $G$ with $d_G(u,v)=\diam(G)$. A diametrical path $P_{d+1}=v_0v_1\cdots v_d$ of a connected graph $G$ with $\diam(G)=d$, is called \textit{separating} if the sets $V(G_0),\dots,V(G_d)$ are pairwise disjoint, that is, if
\begin{equation}\label{eq:decomp}
  V(G)=\bigsqcup_{j=0}^{d}V(G_j),
\end{equation}
where $G_j$ is the connected component of $G-E(P_{d+1})$ containing $v_j$ for $0\le j\le d$. Note that if $P_{d+1}$ is separating, then every path from a vertex of $G_j$ to a vertex outside $G_j$ passes through $v_j$. Moreover, $P_{d+1}$ is separating if and only if every edge of $P_{d+1}$ is a cut edge of $G$. Since every edge of a tree is a cut edge, every diametrical path of a tree is separating.

 The eccentricity matrix of a connected graph $G$, denoted by $\mathcal{E}(G)$, is obtained from the distance matrix of $G$ by keeping each entry that is the largest in its row or the largest in its column, and replacing all other entries by $0$, that is, the $(u,v)$-th entry of $\mathcal{E}(G)$ is defined as $$\mathcal{E}(G)_{uv}=
\begin{cases}
	\text{$d_G(u,v)$} & \quad\text{if $d_G(u,v)=\min\{e_G(u),e_G(v)\}$,}\\
	\text{0} & \quad\text{otherwise.}
\end{cases}$$
Since $\mathcal{E}(G)$ is a real symmetric matrix, all $\mathcal{E}$-eigenvalues of $G$ are real. The eigenvalues of $\mathcal{E}(G)$ are called the $\mathcal{E}$-eigenvalues of $G$, in which the largest one is called the $\mathcal{E}$-spectral radius of $G$. The $\mathcal{E}$-energy of $G$ is the sum of the absolute values of all $\mathcal{E}$-eigenvalues of $G$.

 Recently, the study of eccentricity matrices of graphs has gained significant importance and attracted the attention of researchers.  In \cite{mahato2020spectra}, the authors determined the spectra and inertia of eccentricity matrices of various classes of graphs;  Wei et al. \cite{wei2020solutions} studied the extremal problems for $\mathcal{E}$-spectral radius of trees and determined the trees on $n$ vertices with minimum $\mathcal{E}$-spectral radius; Mahato and Kannan \cite{mahato2022inertia} determined the inertia of eccentricity matrices of trees and characterized the trees for which the $\mathcal{E}$-eigenvalues are symmetric about the origin. For more advances on the eccentricity matrices of graphs, we refer to \cite{lei2022spectral,lei2021eigenvalues,mahato2022eccentricity,mahato2022inertia,mahato2022minimizers,mahato2025characterization,qiu2022eccentricity,wei2022eccentricity} and the references therein.

The eccentricity matrix has numerous applications in chemical graph theory. Randi\'c et al.~\cite{ran2} observed that the eccentricity matrix of a graph (especially a tree) is very sensitive to the branching pattern of molecular graphs, and Wang et al.~\cite{ecc-main} gave some applications of this novel matrix as molecular descriptors. Wang et al.~\cite{wang2020boiling} studied the correlations of the $\mathcal{E}$-spectral radius and the $\mathcal{E}$-energy with the boiling points of benzenoid hydrocarbons and with some chemical properties of octane isomers. So, the study of eccentricity matrix is not only interesting and important from spectral graph theory perspective but also from chemical graph theory perspective.

 An $n \times n$ nonnegative matrix $A $ is  \emph{reducible} if there exists an $n \times n$ permutation matrix $Q$ such that $QAQ^T =
\begin{pmatrix}
A_{11} & A_{12}  \\
0  & A_{22}
\end{pmatrix}$, where $A_{11}$ is an $r \times r$ submatrix with $1 \leq r < n$. If no such permutation matrix $Q$ exists, then $A$ is \emph{irreducible}. Given an $n \times n$ nonnegative symmetric matrix A, associate a graph $G(A)$ on $n$ vertices, say $\{u_1, \dots , u_n\}$ such that there is an edge between the vertices $u_i$ and $u_j$ in $G(A)$ if and only if the $(i,j)^{th}$ entry of $A$ is nonzero. An  $n \times n$  nonnegative symmetric matrix  $A$ is  irreducible if and only if $G(A)$ is connected.

Unlike the adjacency matrix and the distance matrix of a connected graph, the eccentricity matrix need not be irreducible. For example, $\mathcal{E}(K_{a,b})$ is reducible for $a,b\ge2$ while $\mathcal{E}(T)$ is irreducible for every tree $T$ with at least two vertices. So, the following problem is interesting and significant in the context of eccentricity matrix.

\begin{prob}\label{prob:irreducibility}
  Characterize the graphs whose eccentricity matrix is irreducible.    
\end{prob}

Problem~\ref{prob:irreducibility} has been studied for many classes of graphs. Wang et al.~\cite{ecc-main} proved that the eccentricity matrix of every tree is irreducible, and Mahato et al. gave a different proof of this result in~\cite{mahato2020spectra} . In ~\cite{wang2020spectral}, Wang et al. proved the irreducibility of the eccentricity matrix of a connected graph with a cut vertex under a given condition on the cut vertex. Patel et al. \cite{patel2021energy} determined the irreducibility of  windmill graphs, the coalescence of complete graphs, and the coalescence of two cycles.
The eccentricity matrix is irreducible for various classes of graphs, such as clique trees whose blocks contain at most two cut vertices~\cite{li2026inertia_clique_trees}, a class of bi-block graphs that contains all trees with at least three vertices~\cite{divyadevi2025inertia_biblock}, wheels~\cite{jeyaraman2025wheel}, line graphs of starlike trees~\cite{balamoorthy2025starlike}, generalized friendship graphs and generalized corona products~\cite{balamoorthy2025product}.  The irreducibility of the eccentricity matrix has been  determined for various graph operations, such as   $H$-joins~\cite{balamoorthy2025product}, subdivision edge joins~\cite{ashokan2026irreducibility}, central graphs and their joins~\cite{ashokan2026central}, Cartesian products of two trees or of a connected graph and a self-centered graph~\cite{balamoorthy2026irreducibility}, etc. Very recently,  Qiu et al. ~\cite{qiu2026extremal} studied this problem for unicyclic graphs and proved that $\E(U)$ is irreducible if $U$ has a pendant vertex, odd diameter $d\ge3$ and a diametrical path with at most two vertices on the cycle. For more details on the irreducibility of eccentricity matrices of graphs, see the articles  \cite{pandey2025onepoint,pandey2024corona,qiu2023rairo,sorgun2023two} and the references therein.  Despite these results, Problem~\ref{prob:irreducibility} is still open. These works motivate us to study Problem
\ref{prob:irreducibility} for graphs with a separating diametrical path and for unicyclic graphs.

\textit{The inertia} $\In(M)$ of a real symmetric matrix $M$ is the triple $\big(n_{+}(M),n_{-}(M),n_{0}(M)\big)$, where $n_{+}(M),$ $n_{-}(M)$ and $n_{0}(M)$ denote the number of positive, negative and zero eigenvalues of $M$, respectively. There are only a few graphs for which the inertias of eccentricity matrices are known. In \cite{mahato2020spectra}, Mahato et al. computed the inertia of eccentricity matrices of the path and the lollipop graphs; Patel et al. \cite{patel2021energy} studied the inertia of coalescence of graphs. In \cite{mahato2022inertia}, Mahato and Kannan determined the inertia of eccentricity matrices of trees and characterized the trees for which the $\mathcal{E}$-eigenvalues are symmetric about the origin. Recently, in \cite{divyadevi2025inertia_biblock}, Divyadevi and Jeyaraman extended these results to certain bi-block graphs. Very recently, in 2026, Li et al. \cite{li2026inertia_clique_trees} extended these results to clique trees. Motivated by these works, in this article, we consider this problem for the eccentricity matrix of unicyclic graphs and determine  the inertia of the eccentricity matrix of a unicyclic graph whose branches are nontrivial and have equal heights, and prove that its $\mathcal{E}$-spectrum is symmetric about the origin if and only if the girth is even.    

In \cite{mahato2021spectral}, the authors gave a necessary and sufficient condition for the eccentricity matrix of a tree to be invertible. Motivated by this, in this article, we obtain a necessary condition for the invertibility of eccentricity matrix of a unicyclic graph.

The article is organized as follows. In Section~\ref{sec:prelim}, we fix notations and collect some useful results. In Section~\ref{sec:irred}, we give two sufficient conditions for the irreducibility of the eccentricity matrix of a connected graph. Then we show that the eccentricity matrix of every unicyclic graph of girth $3$, $5$ or $7$ is irreducible. We also construct a unicyclic graph of girth $g$ whose eccentricity matrix is reducible for every odd $g\ge9$, and we study the irreducibility of unicyclic graphs whose branches have the same height.  In Section \ref{sec:rank}, we determine the inertia of the eccentricity matrix of a unicyclic graph whose branches are nontrivial and have equal heights, and prove that its $\mathcal{E}$-spectrum is symmetric about the origin if and only if the girth is even.     In Section~\ref{sec:inv}, we obtain a necessary condition for the invertibility of eccentricity matrix of a unicyclic graph.

\section{Notations and preliminaries}\label{sec:prelim}
In this section, we introduce notations and collect some preliminary results which will be used to prove the main results.

 Let ${\mathbb {R}^{n\times n}}$ denote the set of all $n\times n$ matrices with real entries. For $A\in {\mathbb {R}^{n\times n}}$, let  $\det (A)$, $\det(\lambda I-A)$ denote \textit{determinant and characteristic polynomial of $A$}, respectively. Let \textbf{0} denote the \textit{all zero matrix} of appropriate size. For a real number $x$, $\lfloor x \rfloor$ denotes the greatest integer less than or equal to $x$, and $\lceil x \rceil$ denotes the least integer greater than or equal to $x$.

Let $G$ be a simple connected graph. A vertex $w$ is an \emph{eccentric vertex} of $u$ if $d_G(u,w)=e_G(u)$. We omit the subscript $G$ in $d_G$ and $e_G$ when the graph is clear from the context. 
The \emph{eccentric graph} $G_e$ of $G$ is the graph with vertex set $V(G)$ in which two distinct vertices are adjacent if one of them is an eccentric vertex of the other. Since $d(u,w)\le\min\{e(u),e(w)\}$, two distinct vertices $u$ and $w$ are adjacent in $G_e$ if and only if $d(u,w)=\min\{e(u),e(w)\}$, that is, if and only if $\E(G)_{uw}\neq0$. Since $\E(G)$ is a nonnegative symmetric matrix with zero diagonal, and $\E(G)_{uw}\neq0$ exactly when $u$ and $w$ are adjacent in $G_e$, we have the following result.

\begin{thm}[{\cite[Corollary 2.6]{sorgun2023two}}]\label{thm:irr-conn}
Let $G$ be a connected graph on at least two vertices. Then $\E(G)$ is irreducible if and only if $G_e$ is connected.
\end{thm}

The following result is known for the irreducibility of eccentricity matrices of trees.

\begin{thm}[{\cite[Theorem 2.1]{ecc-main}}]\label{thm:tree-irred-known}
The eccentricity matrix $\mathcal{E}(T)$ of a tree with at least two vertices is irreducible.
\end{thm}

For $X\subseteq V(G)$, let $G[X]$ be the subgraph of $G$ induced by $X$, and let $G-X$ be the graph obtained from $G$ by deleting the vertices in $X$. For $v\in V(G)$, let $G-v=G-\{v\}$, and for $F\subseteq E(G)$, let $G-F$ be the graph obtained from $G$ by deleting the edges in $F$. For $X\subseteq V(G)$, let $\chi_X$ be the characteristic vector of $X$, written as a column vector. For $u\in V(G)$, let $R_u$ be the row of $\E(G)$ indexed by $u$. Let $N_G(v)$ be the set of neighbours of $v$ in $G$, and let $\deg_G(v)=|N_G(v)|$ be the degree of $v$. A vertex of degree one is a \emph{pendant vertex}. Let $K_n$, $P_n$, $C_n$ and $K_{1,n-1}$ denote the complete graph, the path, the cycle and the star on $n$ vertices, respectively.

A principal submatrix of a square matrix $A$ is the matrix obtained by deleting any $k$ rows and the corresponding $k$ columns from $A$. If the rows and the columns of $A$ are indexed by $V(G)$ and $X\subseteq V(G)$, then $A[X]$ denotes the principal submatrix of $A$ whose rows and columns are indexed by $X$. The following theorem concerns the interlacing of eigenvalues of a real symmetric matrix and its principal submatrix.  

\begin{thm}[{\cite[Theorem~4.3.28]{hor-john-mat}}, Interlacing Theorem]\label{thm:interlacing}
 Let $A$ be a real symmetric matrix of order $n$, and let $B$ be its principal submatrix of order $m<n$. If $\lambda_1\geq \lambda_2\geq \cdots \geq \lambda_n$ are the eigenvalues of $A$ and $\mu_1\geq \mu_2\geq \cdots \geq \mu_m$ are the eigenvalues of $B$, then $\lambda_{n-m+i} \leq \mu_i \leq \lambda_i$ for $1\leq i \leq m$.
\end{thm}

In the following lemma, we state the well known Schur complement formula.

\begin{lem}[{\cite[Page 4]{bapatbook}}, Schur complement formula]\label{lem:Schur complement}
  Let $A$ be an $n\times n$ matrix partitioned as
\[ A=
\left({\begin{array}{cc}
    A_{11} & A_{12} \\
    A_{21} & A_{22} \\
    \end{array} } \right),\]
where $A_{11}$ and $A_{22}$ are square matrices. If $A_{11}$ is nonsingular, then $\det A= (\det A_{11})\det(A_{22}-A_{21}{A_{11}^{-1}}A_{12})$, and if $A_{22}$ is nonsingular, then  $\det A= (\det A_{22})\det(A_{11}-A_{12}{A_{22}^{-1}}A_{21})$.  
\end{lem}

The following result gives the $\E$-eigenvalues of cycles.

\begin{thm}[{\cite[Theorem~3.2]{ecc-main}}]\label{thm:cycle-spec}
Let $C_n$ be the cycle on $n$ vertices.
\begin{enumerate}
  \item[(i)] If $n=2k$, then the $\E$-eigenvalues of $C_{2k}$ are $k$ and $-k$, each with multiplicity $k$.
  \item[(ii)] If $n=2k+1$, then the $\E$-eigenvalues of $C_{2k+1}$ are $2k\cos\frac{2\pi i}{2k+1}$ for $i=1,\ldots,2k+1$.
\end{enumerate}
\end{thm}

Since $2k+1$ is odd, $\frac{2\pi i}{2k+1}$ is never an odd multiple of $\frac{\pi}{2}$, and so $\cos\frac{2\pi i}{2k+1}\neq0$ for every $i$. Thus, by Theorem~\ref{thm:cycle-spec}, all $\E$-eigenvalues of $C_n$ are nonzero, and $\E(C_n)$ is invertible for every $n\ge3$.

\subsection{Unicyclic graphs}\label{subsec:unicyclic}
Throughout this paper, we denote a unicyclic graph $U$ of order $n$ and girth $g$ by $C_g(T_0,T_1,\ldots,T_{g-1})$, where $C_g=v_0v_1\cdots v_{g-1}v_0$ is the unique cycle of $U$ and $T_i$ is the tree rooted at $v_i$ for each $0\leq i\leq g-1$, that is, 
$T_i$ is obtained as the connected component containing $v_i$ after deleting all edges of $C_g$ from $U$. We call $T_0,T_1,\ldots,T_{g-1}$ the \emph{branches} of $U$. 
The \emph{height} $h_i$ of $T_i$ is the largest distance from $v_i$ to a vertex of $T_i$. The tree $T_i$ is \emph{trivial} if $V(T_i)=\{v_i\}$, that is, if $T_i\cong K_1$. An isomorphism between rooted trees maps the root to the root. For $x\in V(T_i)$, the \emph{children} of $x$ are the neighbours $y$ of $x$ in $T_i$ with $d(y,v_i)=d(x,v_i)+1$, and a \emph{leaf} of $T_i$ is a vertex of $T_i$ other than $v_i$ with no children.\\
Let
\[
  h=\max_{0\le i\le g-1} h_i
  \qquad\text{and}\qquad
  H_i=\max\{d_U(v_i,w) : w\in V(U)\setminus V(T_i)\}.
\]
For each $i$, let $u_i$ be a fixed vertex of $T_i$ with $d_U(u_i,v_i)=h_i$. Note that $u_i=v_i$ when $T_i$ is trivial. The vertex $u_i$ need not be unique, but we use only the equality $d_U(u_i,v_i)=h_i$, and so the results below do not depend on this choice.

Let $\mathcal{U}_{n,g}$ be the set of unicyclic graphs $C_g(T_0,T_1,\ldots,T_{g-1})$ on $n$ vertices in which $T_0,T_1,\ldots,T_{g-1}$ are nontrivial and have the same height. 

Every path from a vertex of $T_i$ to a vertex outside $T_i$ passes through $v_i$. Hence, for $u\in V(T_i)$ and $w\in V(T_j)$ with $i\neq j$,
\begin{equation}\label{eq:branch-dist}
  d(u,w)=d(u,v_i)+d(v_i,v_j)+d(v_j,w).
\end{equation}
In particular, $d(w,u_i)=d(w,v_i)+h_i$ for every $w\notin V(T_i)$. Hence,
\begin{equation}\label{eq:ecc-lower}
 e(w)\ge d(w,v_i)+h_i\qquad\text{for every } w\in V(U)\setminus V(T_i).
\end{equation}

The following lemmas will be used several times in the proofs of subsequent results.

\begin{lem}\label{lem:tree}
Let $U=C_g(T_0,T_1,\ldots,T_{g-1})$ be a unicyclic graph. Let $i\in\{0,1,\ldots,g-1\}$, and assume that $e(v_i)=d(v_i,v)$ for some $v\notin V(T_i)$. Then $e(u)=d(u,v)=d(u,v_i)+e(v_i)$ for every $u\in V(T_i)$. In particular, every vertex of $T_i$ is adjacent to $v$ in $U_e$.
\end{lem}

\begin{proof}
Let $u$ be any vertex in $T_i$. Since every path from $u$ to $v$ passes through $v_i$, we have $$d(u,v)=d(u,v_i)+d(v_i,v)=d(u,v_i)+e(v_i).$$
Therefore, for every $w\in V(U)$,
\[
  d(u,w)\le d(u,v_i)+d(v_i,w)\le d(u,v_i)+e(v_i)=d(u,v),
\]
which implies that $e(u)=d(u,v)$. Thus, $u$ is adjacent to $v$ in $U_e$.
\end{proof}

\begin{lem}\label{lem:struct}
Let $U=C_g(T_0,T_1,\ldots,T_{g-1})$ be a unicyclic graph with $h_0=h_1=\cdots=h_{g-1}=h$, let $k=\lfloor g/2\rfloor$, and let $u\in V(T_i)$. Then $e(u)=d(u,v_i)+h+k$. Moreover, $\E(U)_{uw}=0$ for every $w\in V(T_i)$. For $w\in V(T_j)$ with $i\neq j$,
\begin{align*}
\E(U)_{uw}\ne0
  \quad \text{if and only if} \quad d(v_i,v_j)=k
       \ \text{ and }\ \max\{d(u,v_i),\,d(w,v_j)\}=h,
\end{align*}
and in this case $\E(U)_{uw}=d(u,v_i)+d(w,v_j)+k$.
\end{lem}

\begin{proof}
Observe that for any two vertices $v_i$ and $v_l$ on the cycle $C_g$, we have $d(v_i,v_l)\le k$, and $d(v_i,v_{i+k})=k$ for every $i$. Let $u\in V(T_i)$ and $x\in V(T_l)$. If $l\ne i$, then $d(u,x)=d(u,v_i)+d(v_i,v_l)+d(v_l,x)\le d(u,v_i)+k+h$ by \eqref{eq:branch-dist}, and the equality holds when $d(v_i,v_l)=k$ and $x=u_l$. If $l=i$, then $d(u,x)\le d(u,v_i)+d(v_i,x)\le d(u,v_i)+h$. Hence, $e(u)=d(u,v_i)+h+k$.

Now, let $w\in V(T_j)$, and without loss of generality, assume that $d(u,v_i)\le d(w,v_j)$. Then $e(u)\le e(w)$, and so $\E(U)_{uw}\ne0$ if and only if $d(u,w)=e(u)$. If $j=i$, then $d(u,w)\le d(u,v_i)+h<e(u)$, because $k\ge1$. If $j\ne i$, then, by \eqref{eq:branch-dist}, $d(u,w)=e(u)$ if and only if $d(v_i,v_j)+d(v_j,w)=k+h$. Since $d(v_i,v_j)\le k$ and $d(v_j,w)\le h$, this happens exactly when $d(v_i,v_j)=k$ and $d(w,v_j)=h$, that is, when $d(v_i,v_j)=k$ and $\max\{d(u,v_i),d(w,v_j)\}=h$, because $d(u,v_i)\le d(w,v_j)$. In this case, $\E(U)_{uw}=d(u,w)=d(u,v_i)+k+d(v_j,w)$.
\end{proof}

\section{Irreducibility of the eccentricity matrix}\label{sec:irred}
In this section, we discuss the irreducibility of the eccentricity matrices of graphs with a separating diametrical path and unicyclic graphs, respectively, in the following two subsections.

\subsection{Irreducibility of graphs with a separating diametrical path}
In this subsection, we prove that the eccentricity matrix of a graph containing a separating diametrical path is irreducible. This generalizes Theorem~\ref{thm:tree-irred-known} of \cite{ecc-main}, and the case of \cite[Theorem~2.1]{qiu2026extremal} in which the diametrical path has at most one vertex on the cycle. Furthermore, we give a generalization of Theorem 2.1 of \cite{wang2020spectral}. First, we prove the following lemma to compute the eccentricity of any vertex of a graph containing a  separating diametrical path. This will be used several times in the proof of subsequent results.

\begin{lem}\label{lem:ecc-general}
Let $G$ be a connected graph with $\diam(G)=d$, and let $P_{d+1}=v_0v_1\cdots v_d$ be a separating diametrical path of $G$. Then $e(u)=\max\{d(u,v_0),\,d(u,v_d)\}$ for every $u\in V(G)$.
\end{lem}

\begin{proof}
Without loss of generality, let us assume that  $u\in V(G_i)$. Since $P_{d+1}$ is separating, every path from $u$ to $v_0$ or $u$ to $v_d$ passes through $v_i$. Hence, $d(u,v_0)=d(u,v_i)+d(v_i,v_0)=d(u,v_i)+i$ and $d(u,v_d)=d(u,v_i)+d(v_i,v_d)=d(u,v_i)+(d-i)$. Similarly, for any vertex $w\in V(G_j)$, it follows that $d(w,v_0)=d(w,v_j)+j\le d$ and $d(w,v_d)=d(w,v_j)+(d-j)\le d$. Now, if $j\ge i$, then
\[
  d(u,w)\le d(u,v_i)+d(v_i,v_j)+d(v_j,w)\le d(u,v_i)+(j-i)+(d-j)=d(u,v_d).
\]
Similarly, if $j<i$, then $d(u,w)\le d(u,v_i)+(i-j)+j=d(u,v_0)$. Thus, $e(u)\le\max\{d(u,v_0),d(u,v_d)\}\le e(u)$.
\end{proof}

In the following theorem, we prove that the eccentricity matrix of a graph containing a separating diametrical path is irreducible. 

\begin{thm}\label{thm:branch-irred}
Let $G$ be a connected graph with $\diam(G)=d\ge1$, and assume that $G$ has a separating diametrical path $P_{d+1}=v_0v_1\cdots v_d$. Then $\mathcal{E}(G)$ is irreducible.
\end{thm}
\begin{proof}
For any vertex $u\in V(G)$, by Lemma~\ref{lem:ecc-general}, it follows that $e(u)=\max\{d(u,v_0),\,d(u,v_d)\}$. Therefore, every vertex of $G$ is adjacent to at least one of $v_0$ or $v_d$ in $G_e$, and hence $G_e$ is connected. Thus, $\mathcal{E}(G)$ is irreducible by Theorem~\ref{thm:irr-conn}.
\end{proof}

\begin{rmk}\label{rmk:tree-irred}
Observe that every diametrical path of a tree is separating. Therefore, Theorem \ref{thm:branch-irred} generalizes Theorem \ref{thm:tree-irred-known} of \cite{ecc-main}.
\end{rmk}

Qiu et al.~\cite[Theorem~2.1]{qiu2026extremal} proved that if a unicyclic graph $U$ has at least one pendant vertex, odd diameter $d\ge3$ and a diametrical path with at most two vertices on $C_g$, then $\E(U)$ is irreducible. Corollary~\ref{cor:unicyclic-branch-irred} has no condition on the diameter, but it requires a diametrical path with at most one vertex on $C_g$.

\begin{cor}\label{cor:unicyclic-branch-irred} 
Let $U$ be a unicyclic graph. If the cycle $C_g$ of $U$ contains at most one vertex of some diametrical path of $U$, then $\mathcal{E}(U)$ is irreducible.
\end{cor}

\begin{proof}
Let $P_{d+1}$ be a diametrical path of $U$ with at most one vertex on $C_g$. Then no edge of $P_{d+1}$ is an edge of $C_g$. Since $C_g$ is the only cycle of $U$, every edge of $U$ not on $C_g$ is a cut edge. Thus, every edge of $P_{d+1}$ is a cut edge of $U$, and $P_{d+1}$ is separating. The result now follows from Theorem~\ref{thm:branch-irred}.
\end{proof}

\begin{thm}\label{thm:cutvertex_irreducible}
Let $G$ be a connected graph with a cut vertex $c$, and let $B$ be a connected component of $G-c$ such that $\hat B=G[V(B)\cup\{c\}]$ has a separating diametrical path. If $e(c)=d(c,v)$ for some $v\in V(B)$, then $\mathcal{E}(G)$ is irreducible.
\end{thm}

\begin{proof}
Let $X=V(G)\setminus V(B)$. Then $c\in X$. Since $B$ is a connected component of $G-c$, every path from a vertex of $B$ to a vertex of $X$ passes through $c$. Thus, every shortest path between two vertices of $\hat B$ lies in $\hat B$, and so $d_{\hat B}(u,w)=d(u,w)$ for all $u,w\in V(\hat B)$. Since $e(c)=d(c,v)$ with $v\in V(B)$, we have $e_{\hat B}(c)=e(c)$. Let $a$ and $b$ be the end vertices of a separating diametrical path of $\hat B$. By Lemma~\ref{lem:ecc-general}, $e_{\hat B}(u)=\max\{d(u,a),d(u,b)\}$ for every $u\in V(\hat B)$. Without loss of generality, assume that $d(c,b)=e_{\hat B}(c)=e(c)$. Then $b\neq c$, and so $b\in V(B)$.

By Theorem~2.1, it suffices to show that $G_e$ is connected, which we establish through the following three claims.

\smallskip
\noindent\textbf{Claim 1. Every $x\in X$ is adjacent to $b$ in $G_e$.}

\noindent\emph{Proof of Claim~1.} Let $x\in X$. Since every path from $x$ to $b$ passes through $c$, we have $d(x,b)=d(x,c)+d(c,b)=d(x,c)+e(c)$. On the other hand, $d(x,w)\le d(x,c)+d(c,w)\le d(x,c)+e(c)$ for every $w\in V(G)$, and so $e(x)\le d(x,c)+e(c)=d(x,b)$. Thus, $d(x,b)=e(x)$, that is, $b$ is an eccentric vertex of $x$. Since $x\neq b$, the vertices $x$ and $b$ are adjacent in $G_e$.

\smallskip
\noindent\textbf{Claim 2. Every $u\in V(B)$ is adjacent in $G_e$ to a vertex of $X\cup\{a,b\}$.}

\noindent\emph{Proof of Claim~2.} Let $u\in V(B)$, and let $w$ be an eccentric vertex of $u$. If $w\in X$, then $u$ is adjacent to $w$ in $G_e$. Otherwise $w\in V(B)$. Since $u,w\in V(\hat B)$, we get $e(u)=d(u,w)\le e_{\hat B}(u)\le e(u)$. Thus, $e(u)=e_{\hat B}(u)=\max\{d(u,a),d(u,b)\}$, and $a$ or $b$ is an eccentric vertex of $u$. Since $e(u)\ge1$, this eccentric vertex is different from $u$, and so $u$ is adjacent to it in $G_e$.

\smallskip
\noindent\textbf{Claim 3. The vertices $a$ and $b$ lie in the same connected component of $G_e$.}

\noindent\emph{Proof of Claim~3.} If $a=c$, then $a\in X$, and $a$ is adjacent to $b$ by Claim~1. Now, let $a\neq c$. Then $a\in V(B)$, and by Claim~2, $a$ is adjacent in $G_e$ to a vertex $y\in X\cup\{a,b\}$. Since $y\neq a$, either $y=b$, or $y\in X$ and $y$ is adjacent to $b$ by Claim~1. In both cases, $a$ and $b$ lie in the same connected component of $G_e$.

\smallskip
By Claims~1 and~3, all vertices of $X\cup\{a,b\}$ lie in the connected component of $G_e$ containing $b$. By Claim~2, every vertex of $B$ also lies in this component. Since $V(G)=X\cup V(B)$, the graph $G_e$ is connected, and hence $\mathcal{E}(G)$ is irreducible by Theorem~\ref{thm:irr-conn}.
\end{proof}

From Theorem \ref{thm:cutvertex_irreducible}, we get Theorem 2.1 of \cite{wang2020spectral} as an immediate corollary. Note that there is an ambiguity in the statement of Theorem 2.1 of \cite{wang2020spectral}, which we correct in the following corollary.

\begin{cor}\label{cor:cutvertex-tree}
Let $G$ be a connected graph with a cut vertex $c$, and let $T$ be a connected component of $G-c$ such that $\hat T=G[V(T)\cup\{c\}]$ is a tree. If $e(c)=d(c,v)$ for some $v\in V(T)$, then $\mathcal{E}(G)$ is irreducible.
\end{cor}

\begin{proof}
Since every diametrical path of the tree $\hat T$ is separating, the result follows from Theorem~\ref{thm:cutvertex_irreducible}.
\end{proof}

\subsection{Irreducibility of unicyclic graphs}\label{subsec:odd-girth}
In this subsection, we show that the eccentricity matrix of every unicyclic graph of girth $3$, $5$ or $7$ is irreducible but this does not hold for any larger odd girth. For every odd girth $g\ge 9$, we construct a unicyclic graph whose eccentricity matrix is reducible. On the other hand, we also provide a class of unicyclic graphs of odd girth, whose eccentricity matrices are always irreducible.  First, we prove the following lemma about the structure of eccentric graph of a unicyclic graph of odd girth.

\begin{lem}\label{lem:odd-facts}
Let $U=C_g(T_0,T_1,\ldots,T_{g-1})$, where $g=2k+1$ is odd, and assume that $e(v_i)>h_i$ for every $i$. Then the following hold.
\begin{enumerate}
  \item[(i)] For every $i$, there is an index $j\neq i$ such that every vertex of $T_i$ is adjacent to $u_j$ in $U_e$.
  \item[(ii)] If $h_j=h$ and $d(v_i,v_j)=k$, then every vertex of $T_i$ is adjacent to $u_j$ in $U_e$.
  \item[(iii)] If $h_j=h_{j+1}=h$, then $u_j$ and $u_{j+1}$ have a common neighbour in $U_e$.
  \item[(iv)] If there is an index $j\neq i$ with $h_j=h$ and $d(v_i,v_j)\ge k-1$, then there is an index $l\neq i$ with $h_l=h$ such that every vertex of $T_i$ is adjacent to $u_l$ in $U_e$.
\end{enumerate}
\end{lem}

\begin{proof}
Since $e(v_i)>h_i$, we have
\[
e(v_i)=H_i=\max_{j\neq i}d(v_i,u_j) \quad\text{for every } i.
\]

(i) Taking $j\neq i$ with $e(v_i)=d(v_i,u_j)$, we get the result from Lemma~\ref{lem:tree}.

(ii) For every $l\neq i$, we have $d(v_i,u_l)=d(v_i,v_l)+h_l\le k+h=d(v_i,u_j)$. Thus, $e(v_i)=d(v_i,u_j)$, and Lemma~\ref{lem:tree} gives the result.

(iii) We have $d(v_{j+k+1},v_j)=d(v_{j+k+1},v_{j+1})=k$. By (ii) with $i=j+k+1$, the vertex $v_{j+k+1}$ is adjacent to both $u_j$ and $u_{j+1}$ in $U_e$.

(iv) If $l\neq i$ and $h_l<h$, then $d(v_i,u_l)\le k+h-1\le d(v_i,v_j)+h=d(v_i,u_j)$. Thus, the maximum of $d(v_i,u_l)$ over all $l\neq i$ is attained at some $l$ with $h_l=h$. For this $l$, we have $e(v_i)=d(v_i,u_l)$, and Lemma~\ref{lem:tree} gives the result.
\end{proof}
In the following theorem, we prove that the eccentricity matrix of any unicyclic graph with girth $3,5$ or $7$ is irreducible.

\begin{thm}\label{thm:girth357}
Let $U$ be a unicyclic graph with girth $g$. If $g\in\{3,5,7\}$, then $\E(U)$ is irreducible.
\end{thm}

\begin{proof}
By Theorem~\ref{thm:irr-conn}, it is enough to show that $U_e$ is connected. Let $U=C_g(T_0,T_1,\ldots,T_{g-1})$ with $g=2k+1$, where $k\in\{1,2,3\}$, and let $J=\{j: h_j=h\}$. Without loss of generality, assume that $h_0=h$. We consider two cases. If $e(v_0)=h$, then the result follows from Corollary~\ref{cor:cutvertex-tree}. If $e(v_0)>h$, then we show that the connected component of $U_e$ containing $u_0$ contains every vertex of $U$. For this, we prove four claims. In Claim~1, we show that this component contains $u_j$ for every $j\in J$. In Claims~2 and~3, we show that it contains $V(T_i)$ for every $i\neq0$. In Claim~4, we show that it contains $V(T_0)$.

\smallskip
\noindent\textbf{Case 1.} Let $e(v_0)=h$. Since $e(v_0)\ge 1$, we have $h\ge1$, and the tree $T_0$ is nontrivial. Thus, $v_0$ is a cut vertex of $U$. Let $T$ be the connected component of $U-v_0$ containing $u_0$. Then $U[V(T)\cup\{v_0\}]$ is a subtree of $T_0$, and $e(v_0)=h=d(v_0,u_0)$. Therefore, by Corollary~\ref{cor:cutvertex-tree}, $\E(U)$ is irreducible.

\smallskip
\noindent\textbf{Case 2.} Let $e(v_0)>h$. For $i\neq 0$, we have $e(v_i)\ge d(v_i,u_0)=d(v_i,v_0)+h>h\ge h_i$. Thus, $e(v_i)>h_i$ for every $i$, and Lemma~\ref{lem:odd-facts} applies. Let $X$ be the vertex set of the connected component of $U_e$ containing $u_0$. We prove the following four claims.

\smallskip
\noindent\textbf{Claim 1. $u_j\in X$ for every $j\in J$.}

\noindent\emph{Proof of Claim~1.} If $|J|=1$, then $J=\{0\}$, and Claim~1 holds because $u_0\in X$. Thus, assume that $|J|\ge2$. For distinct $j,j'\in J$, the vertices $u_j$ and $u_{j'}$ lie in the same connected component of $U_e$ whenever $d(v_j,v_{j'})\in\{1,k\}$. Indeed, if $d(v_j,v_{j'})=k$, then $u_j$ and $u_{j'}$ are adjacent in $U_e$ by Lemma~\ref{lem:odd-facts}(ii), and if $d(v_j,v_{j'})=1$, then they have a common neighbour in $U_e$ by Lemma~\ref{lem:odd-facts}(iii). For $g\in\{3,5\}$, any two distinct vertices of $C_g$ are at distance $1$ or $k$. Since $0\in J$ and $u_0\in X$, Claim~1 follows in this case.

Now, let $g=7$, so that $k=3$, and suppose that Claim~1 fails. Let $J_1=\{j\in J: u_j\in X\}$ and $J_2=J\setminus J_1$. Then $0\in J_1$ and $J_2\neq\emptyset$. By the previous paragraph, $d(v_p,v_q)=2$ for all $p\in J_1$ and $q\in J_2$. Taking $q\in J_2$, we get $J_1\subseteq\{q-2,q+2\}$, since $v_{q-2}$ and $v_{q+2}$ are the only vertices of $C_7$ at distance $2$ from $v_q$. We claim that $J_1$ or $J_2$ has exactly one element. If $|J_1|=1$, there is nothing to prove. Otherwise $J_1=\{q-2,q+2\}$, and every $q'\in J_2$ is at distance $2$ from both $v_{q-2}$ and $v_{q+2}$. The only such vertex of $C_7$ is $v_q$, so $J_2=\{q\}$.

Let $j$ be the element of the part that has exactly one element. Since both $J_1$ and $J_2$ are nonempty, the set $J\setminus\{j\}$ is nonempty, and every $l\in J\setminus\{j\}$ lies in the other part. Thus, $d(v_j,v_l)=2=k-1$ for every $l\in J\setminus\{j\}$. By Lemma~\ref{lem:odd-facts}(iv) with $i=j$, there is an index $l\in J\setminus\{j\}$ such that every vertex of $T_j$, in particular $u_j$, is adjacent to $u_l$ in $U_e$. Since $j$ and $l$ lie in different parts, exactly one of $u_j$ and $u_l$ lies in $X$. This contradicts the fact that $X$ is the vertex set of a connected component of $U_e$. Hence Claim~1 holds.

\smallskip
\noindent\textbf{Claim 2. If $i\neq0$ and $d(v_i,v_0)\ge k-1$, then $V(T_i)\subseteq X$.}

\noindent\emph{Proof of Claim~2.} By Lemma~\ref{lem:odd-facts}(iv) with $j=0$, every vertex of $T_i$ is adjacent in $U_e$ to $u_l$ for some $l\in J$, and $u_l\in X$ by Claim~1. For $g\in\{3,5\}$, we have $k-1\le1$, so Claim~2 covers every $i\neq0$.

\smallskip
\noindent\textbf{Claim 3. If $g=7$ and $i\in\{1,6\}$, then $V(T_i)\subseteq X$.}

\noindent\emph{Proof of Claim~3.} Without loss of generality, assume that $i=1$. We have $e(v_1)\ge d(v_1,u_0)=h+1$. If $e(v_1)=h+1$, then $e(v_1)=d(v_1,u_0)$, and every vertex of $T_1$ is adjacent to $u_0$ by Lemma~\ref{lem:tree}. Otherwise $e(v_1)=d(v_1,u_l)\ge h+2$ for some $l\neq1$. Since $h_l\le h$, we get $d(v_1,v_l)\ge2$, that is, $l\in\{3,4,5,6\}$. If $l=6$, then $d(v_1,v_6)=2$ gives $h_6=h$, and $u_6\in X$ by Claim~1. If $l\in\{3,4,5\}$, then $d(v_l,v_0)\ge2$, and $u_l\in X$ by Claim~2. In both cases, every vertex of $T_1$ is adjacent to $u_l\in X$ by Lemma~\ref{lem:tree}.

\smallskip
\noindent\textbf{Claim 4. $V(T_0)\subseteq X$.}

\noindent\emph{Proof of Claim~4.} By Lemma~\ref{lem:odd-facts}(i), every vertex of $T_0$ is adjacent in $U_e$ to $u_j$ for some $j\neq0$, and $u_j\in X$ by Claims~2 and~3.

\smallskip
Using Claims~1--4, we get $V(T_i)\subseteq X$ for every $i$. Thus, $V(U)=\bigcup_{i}V(T_i)\subseteq X$, the graph $U_e$ is connected, and $\E(U)$ is irreducible.
\end{proof}

Note that Theorem~\ref{thm:girth357} cannot be extended to a larger odd girth. In the following theorem, we give a construction of unicyclic graphs of girth $g$, whose eccentricity matrix is reducible for every odd $g\ge9$.

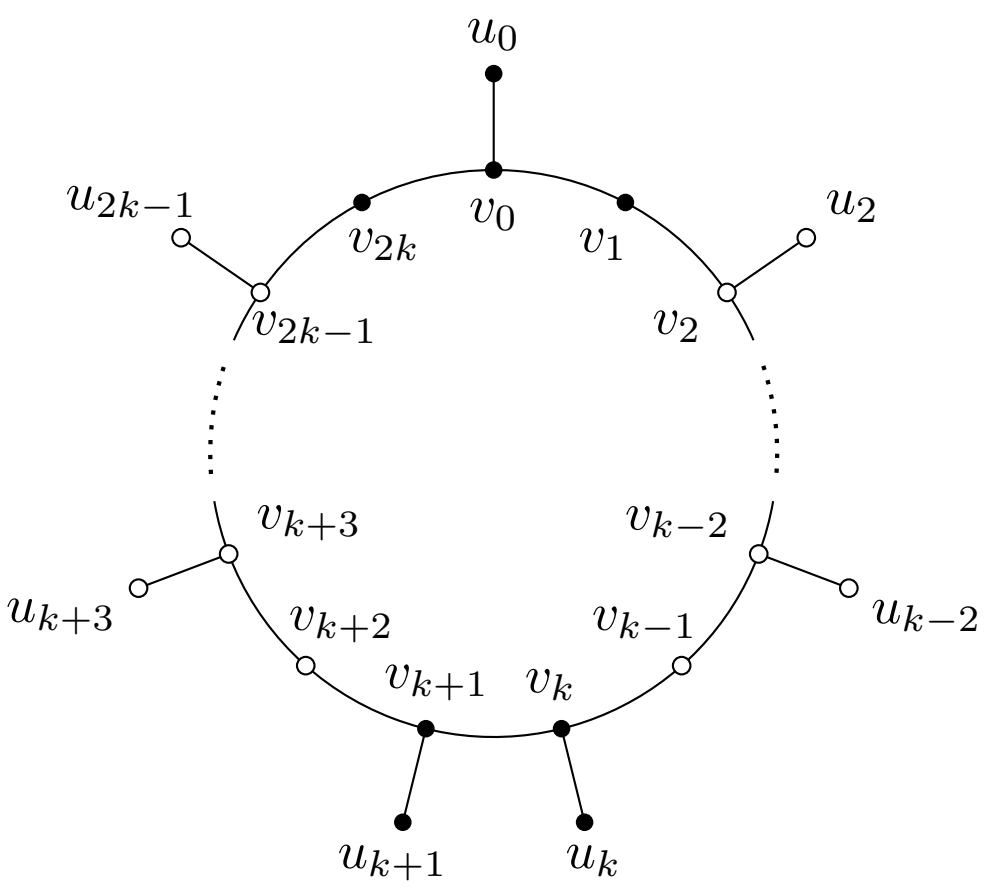
\begin{figure}[htbp]
\centering
\begin{tikzpicture}[scale=0.75]
  \def\r{2.5}\def\R{3.35}
  \pgfmathsetmacro{\st}{360/13}
  \draw (90+2*\st:\r) arc[start angle=90+2*\st, end angle=90-2*\st, radius=\r];
  \draw (90-4*\st:\r) arc[start angle=90-4*\st, end angle=90-9*\st, radius=\r];
  \foreach \a in {2,9} {
    \draw (90-\a*\st:\r) arc[start angle=90-\a*\st, end angle=90-\a*\st-0.4*\st, radius=\r];
    \draw[dotted, thick] (90-\a*\st-0.6*\st:\r) arc[start angle=90-\a*\st-0.6*\st, end angle=90-\a*\st-1.4*\st, radius=\r];
    \draw (90-\a*\st-1.6*\st:\r) arc[start angle=90-\a*\st-1.6*\st, end angle=90-\a*\st-2*\st, radius=\r];
  }
  \foreach \n in {0,2,4,6,7,9,11} {\draw ({90-\n*\st}:\r)--({90-\n*\st}:\R);}
  \foreach \n in {0,1,6,7,12} {\fill ({90-\n*\st}:\r) circle (2.2pt);}
  \foreach \n in {0,6,7} {\fill ({90-\n*\st}:\R) circle (2.2pt);}
  \foreach \n in {2,4,5,8,9,11} {\filldraw[fill=white] ({90-\n*\st}:\r) circle (2.2pt);}
  \foreach \n in {2,4,9,11} {\filldraw[fill=white] ({90-\n*\st}:\R) circle (2.2pt);}
  \foreach \n/\l in {0/0,1/1,2/2,4/k-2,5/k-1,6/k,7/k+1,8/k+2,9/k+3,11/2k-1,12/2k}
    {\node[anchor={90-\n*\st}, inner sep=1pt] at ({90-\n*\st}:\r-0.2) {\small $v_{\l}$};}
  \foreach \n/\l in {0/0,2/2,4/k-2,6/k,7/k+1,9/k+3,11/2k-1}
    {\node[anchor={270-\n*\st}, inner sep=1pt] at ({90-\n*\st}:\R+0.15) {\small $u_{\l}$};}
\end{tikzpicture}
\caption{The graph $U$ of Theorem~\ref{thm:odd-reducible}, of order $4k-2$. The black vertices form the set $W$ of the proof. No edge of $U_e$ joins a black vertex and a white vertex.}
\label{fig:odd-reducible}
\end{figure}

\begin{thm}\label{thm:odd-reducible}
Let $g=2k+1\ge9$, and let $U=C_g(T_0,T_1,\ldots,T_{g-1})$, where $T_i\cong K_1$ for $i\in\{1,k-1,k+2,2k\}$ and $T_i\cong P_2$ for all other $i$ (see Figure~\ref{fig:odd-reducible}). Then $\E(U)$ is reducible.
\end{thm}

\begin{proof}
By Theorem~\ref{thm:irr-conn}, it is enough to show that $U_e$ is disconnected. Let $A=\{0,1,k,k+1,2k\}$ and $W=\bigcup_{i\in A}V(T_i)$. We show that no edge of $U_e$ joins a vertex of $W$ to a vertex outside $W$. If we show this, then we are done, as $v_0\in W$ but $v_2\notin W$. To show this, we first find the eccentricity of every vertex of $U$, and then we find the edges of $U_e$.

On the cycle $C_g$, we have $d(v_i,v_j)\le k$ for all $i,j$, and $d(v_i,v_j)=k$ if and only if $j\in\{i+k,i+k+1\}$. Since $k\ge4$, no two of $v_1,v_{k-1},v_{k+2},v_{2k}$ are adjacent. As $v_{i+k}$ and $v_{i+k+1}$ are adjacent, for every $i$ there is an index $j$ with $d(v_i,v_j)=k$ and $T_j\cong P_2$.

Let $u\in V(T_i)$. For $w\in V(T_i)$, we have $d(u,w)\le1$. For $w\in V(T_j)$ with $j\neq i$, \eqref{eq:branch-dist} gives $d(u,w)\le d(u,v_i)+k+1$, with equality when $d(v_i,v_j)=k$ and $w$ is the pendant vertex of $T_j\cong P_2$. Thus, $e(u)=d(u,v_i)+k+1$ for every $u\in V(T_i)$, as for every $i$ there is an index $j$ with $d(v_i,v_j)=k$ and $T_j\cong P_2$.

Now, let $u\in V(T_i)$ and $w\in V(T_j)$ be adjacent in $U_e$, and without loss of generality, assume that $d(u,v_i)\le d(w,v_j)$. Then $e(u)\le e(w)$, and so $d(u,w)=e(u)$. If $i=j$, then $d(u,w)\le1<e(u)$, which is impossible. Thus, $i\neq j$, and \eqref{eq:branch-dist} gives
\[
  d(u,v_i)+d(v_i,v_j)+d(v_j,w)=d(u,w)=e(u)=d(u,v_i)+k+1,
\]
that is, $d(v_i,v_j)+d(v_j,w)=k+1$. Since $d(v_i,v_j)\le k$ and $d(v_j,w)\le1$, we get $d(v_i,v_j)=k$ and $T_j\cong P_2$. Hence every edge of $U_e$ joins a vertex of $T_i$ and a vertex of $T_j$, where $d(v_i,v_j)=k$ and at least one of $T_i$ and $T_j$ is isomorphic to $P_2$.

Finally, for $i=0,1,k,k+1,2k$, the vertices at distance $k$ from $v_i$ are $\{v_k,v_{k+1}\}$, $\{v_{k+1},v_{k+2}\}$, $\{v_{2k},v_0\}$, $\{v_0,v_1\}$ and $\{v_{k-1},v_k\}$, respectively. Since $k\ge4$, we have $k-1\notin A$ and $k+2\notin A$. Thus, the only pairs $\{i,j\}$ with $i\in A$, $j\notin A$ and $d(v_i,v_j)=k$ are $\{1,k+2\}$ and $\{2k,k-1\}$. Since $T_i\cong K_1$ for $i\in\{1,k-1,k+2,2k\}$, no edge of $U_e$ joins a vertex of $W$ to a vertex outside $W$, and $\E(U)$ is reducible.
\end{proof}

\begin{cor}\label{cor:girth-char}
Let $g\ge3$. Then $\E(U)$ is irreducible for every unicyclic graph $U$ of girth $g$ if and only if $g\in\{3,5,7\}$.
\end{cor}
 In the following theorem, we provide a class of unicyclic graphs of odd girth, whose eccentricity matrices are always irreducible.
 
\begin{thm}\label{thm:oddgirth-equal}
Let $U=C_g(T_0,T_1,\ldots,T_{g-1})$ be a unicyclic graph with girth $g=2k+1$ such that all $T_0,T_1,\ldots,T_{g-1}$ have equal height. Then $\mathcal{E}(U)$ is irreducible.
\end{thm}

\begin{proof}
By Theorem~\ref{thm:irr-conn}, it is enough to show that $U_e$ is connected. We first show that every vertex of $T_i$ is adjacent to $u_{i+k}$ in $U_e$, and then we show that the vertices $u_0,u_1,\ldots,u_{g-1}$ lie in one connected component of $U_e$.

Let $u\in V(T_i)$. Since $d(v_i,v_{i+k})=k$ and $d(u_{i+k},v_{i+k})=h$, Lemma~\ref{lem:struct} gives $u$ is adjacent to $u_{i+k}$ in $U_e$. In particular, $u_i$ is adjacent to $u_{i+k}$ for every $i$. Since $\gcd(k,2k+1)=1$, the walk $u_0,u_k,u_{2k},\ldots,u_{(g-1)k}$ in $U_e$ passes through every $u_j$. Thus, all the vertices $u_j$ lie in one connected component of $U_e$, and every other vertex is adjacent to one of them. Hence $U_e$ is connected, and $\mathcal{E}(U)$ is irreducible.
\end{proof}


 Note that the result of Theorem \ref{thm:oddgirth-equal} does not hold for even girth. In the following theorem, we prove that the eccentricity matrix of any unicyclic graph $U=C_g(T_0,T_1,\ldots,T_{g-1})$ of even girth, where $T_0,T_1,\ldots,T_{g-1}$ are of equal height, is reducible. 

\begin{thm}\label{thm:evengirth_reducible}
Let $U=C_g(T_0,T_1,\ldots,T_{g-1})$, where $g=2k\ge4$ is even and $T_0,T_1,\ldots,T_{g-1}$ have the same height. Then $\mathcal{E}(U)$ is reducible.
\end{thm}

\begin{proof}
By Theorem~\ref{thm:irr-conn}, it suffices to show that $U_e$ is disconnected. By Lemma~\ref{lem:struct}, if $u\in V(T_i)$ and $w\in V(T_j)$ are adjacent in $U_e$, then $d(v_i,v_j)=k$, that is, $j=i+k$. Thus, for every $i$, no edge of $U_e$ joins $V(T_i)\cup V(T_{i+k})$ to a vertex outside this set, and $U_e$ has at least $k\ge2$ connected components. Thus, $U_e$ is disconnected and hence $\mathcal{E}(U)$ is reducible.
\end{proof}

Combining Theorems~\ref{thm:oddgirth-equal} and~\ref{thm:evengirth_reducible}, we get the following corollary.

\begin{cor}\label{cor:dichotomy}
Let $U=C_g(T_0,T_1,\ldots,T_{g-1})$, where $T_0,T_1,\ldots,T_{g-1}$ have the same height. Then $\mathcal{E}(U)$ is irreducible if and only if $g$ is odd. In particular, $\E(C_g)$ is irreducible if and only if $g$ is odd.
\end{cor}

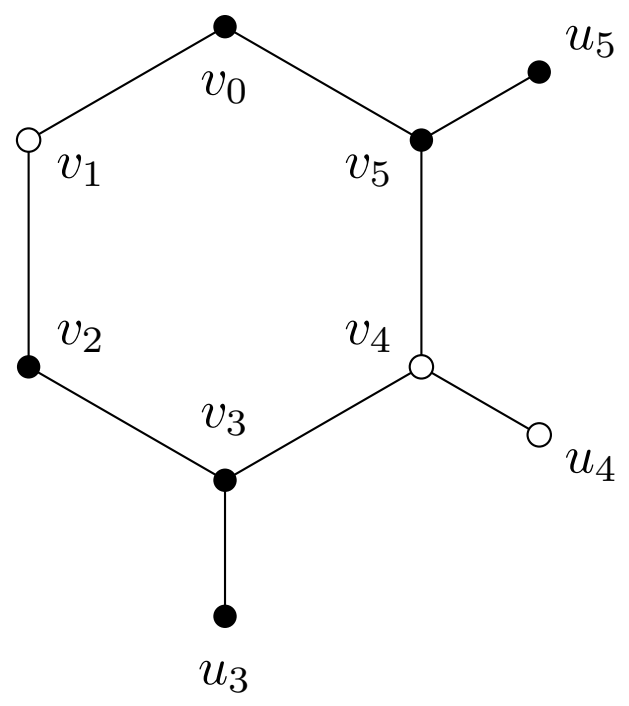
\begin{figure}[htbp]
\centering
\begin{tikzpicture}[scale=1]
  \foreach \i in {0,...,5} {\coordinate (v\i) at ({90+60*\i}:1.5);}
  \foreach \i in {3,4,5} {\coordinate (p\i) at ({90+60*\i}:2.4);}
  \draw (v0)--(v1)--(v2)--(v3)--(v4)--(v5)--cycle;
  \foreach \i in {3,4,5} {\draw (v\i)--(p\i);}
  \foreach \i in {0,2,3,5} {\fill (v\i) circle (2.2pt);}
  \foreach \i in {3,5} {\fill (p\i) circle (2.2pt);}
  \foreach \i in {1,4} {\filldraw[fill=white] (v\i) circle (2.2pt);}
  \filldraw[fill=white] (p4) circle (2.2pt);
  \foreach \i in {0,...,5} {\node at ({90+60*\i}:1.1) {\small $v_{\i}$};}
  \foreach \i in {3,4,5} {\node at ({90+60*\i}:2.8) {\small $u_{\i}$};}
\end{tikzpicture}
\caption{The graph $U=C_6(K_1,K_1,K_1,P_2,P_2,P_2)$}
\label{fig:counter2_graph}
\end{figure}

\begin{rmk}
The converse of Theorem~\ref{thm:evengirth_reducible} does not hold, that is, a unicyclic graph of even girth with a reducible eccentricity matrix need not have branches of the same height. For example, let $U=C_6(K_1,K_1,K_1,P_2,P_2,P_2)$ (see Figure~\ref{fig:counter2_graph}). Not all branches of $U$ have the same height. Here $e(v_3)=e(v_4)=e(v_5)=3$, and every other vertex of $U$ has eccentricity $4$. Thus, $U_e$ has two connected components, with vertex sets $\{v_0,v_2,v_3,v_5,u_3,u_5\}$ and $\{v_1,v_4,u_4\}$. Hence, $\E(U)$ is reducible.
\end{rmk}

\section{Rank and inertia}\label{sec:rank}

In this section, we find the rank and inertia of eccentricity matrices of unicyclic graphs $U=C_g(T_0,T_1,\ldots,T_{g-1})\in\mathcal{U}_{n,g}$, whose branches are nontrivial and have the same height $h\ge1$. Moreover, we characterize the unicyclic graphs in $\mathcal{U}_{n,g}$, whose $\E$-eigenvalues are symmetric with respect to the origin. First, we define a matrix $M$ that will be used in the proof of the next theorem.

Let $M=A\big((C_g)_e\big)$ be the adjacency matrix of the eccentric graph of $C_g$. Since every vertex of $C_g$ has eccentricity $k$, we have $\E(C_g)=kM$. From Theorem \ref{thm:cycle-spec}, it follows that the $\E$-eigenvalues of $C_g$ are nonzero and hence $\E(C_g)$ is invertible. Consequently, $M$ is also invertible.


\begin{thm}\label{thm:rank-inertia}
Let $U=C_g(T_0,T_1,\ldots,T_{g-1})\in\mathcal{U}_{n,g}$. Then $\rank\E(U)=2g$ and $\In\E(U)=(g,g,n-2g)$.
\end{thm}

\begin{proof}
Let $\E=\E(U)$, and let $Z=\{v_0,\ldots,v_{g-1},u_0,\ldots,u_{g-1}\}$. Since $h\ge1$, the set $Z$ has $2g$ vertices. We first show that every row of $\E$ is a multiple of one of the $2g$ rows $R_{v_0},\ldots,R_{v_{g-1}},R_{u_0},\ldots,R_{u_{g-1}}$. 

Let $u\in V(T_i)$ and $w\in V(T_j)$. By Lemma~\ref{lem:struct}, $\E_{uw}\ne0$ if and only if $d(v_i,v_j)=k$ and $\max\{d(u,v_i),d(w,v_j)\}=h$, and in this case $\E_{uw}=d(u,v_i)+d(w,v_j)+k$. If $d(u,v_i)=h$, then $\E_{uw}\ne0$ if and only if $d(v_i,v_j)=k$, and then $\E_{uw}=h+d(w,v_j)+k$. Thus, $R_u=R_{u_i}$. If $d(u,v_i)<h$, then $R_u=(d(u,v_i)+h+k)\chi_Y^T$, where $Y$ is the set of vertices $w$ that lie in a branch $T_j$ with $d(v_i,v_j)=k$ and satisfy $d(w,v_j)=h$. In particular, $R_{v_i}=(h+k)\chi_Y^T$, and so $R_u$ is a multiple of $R_{v_i}$. Therefore, $\rank\E\le2g$.

Now, we show that $\E[Z]$ is nonsingular with inertia $(g,g,0)$. Since $d(v_i,v_i)=0<h$ and $d(u_i,v_i)=h$, we also get $\E_{v_iv_j}=0$, $\E_{v_iu_j}=(h+k)M_{ij}$ and $\E_{u_iu_j}=(2h+k)M_{ij}$, where $M_{ij}$ is the $(i,j)$-entry of $M$. Thus, with the vertices of $Z$ in the order $v_0,\ldots,v_{g-1},u_0,\ldots,u_{g-1}$,
\[
  \E[Z]=\begin{pmatrix}\textbf{0}&(h+k)M\\[2pt](h+k)M&(2h+k)M\end{pmatrix}.
\]
Since $M$ is nonsingular, by Schur complement formula, we have
\begin{align*}
  \det\E[Z]&= \det ((2h+k)M) \det[\textbf{0}-(h+k)M(2h+k)^{-1}M^{-1}(h+k)M]\\
  &=(-1)^{g}(h+k)^{2g}(\det M)^{2}  
\end{align*}

Since $\det M\neq 0$, $\det\E[Z]\neq 0$ and hence $\rank\E=2g$. The principal submatrix of $\E[Z]$ indexed by $v_0,\ldots,v_{g-1}$ is the zero matrix of order $g$. Thus, Theorem~\ref{thm:interlacing}, with $i=1$ and $i=g$, gives $\lambda_{g+1}(\E[Z])\le0\le\lambda_g(\E[Z])$. Since $\E[Z]$ is nonsingular, both inequalities are strict, and so $\In\E[Z]=(g,g,0)$. If $n=2g$, then $\E=\E[Z]$, and the result follows. Otherwise, Theorem~\ref{thm:interlacing} gives $\lambda_g(\E)\ge\lambda_g(\E[Z])>0$ and $\lambda_{n-g+1}(\E)\le\lambda_{g+1}(\E[Z])<0$. Hence $n_+(\E)\ge g$ and $n_-(\E)\ge g$. Since $n_+(\E)+n_-(\E)=2g$, we get $\In\E=(g,g,n-2g)$.
\end{proof}

Next, we characterize the unicyclic graphs in $\mathcal{U}_{n,g}$, whose $\E$-eigenvalues of $U$ are symmetric about the origin. 

\begin{thm}\label{thm:symmetry}
Let $U\in\mathcal{U}_{n,g}$. Then the $\E$-eigenvalues of $U$ are symmetric about the origin if and only if $g$ is even.
\end{thm}

\begin{proof}
Let $\E=\E(U)$. We prove this theorem by considering two cases. For even $g$, we show that $\E$ is similar to $-\E$. For odd $g$, we show that $\operatorname{tr}\E^{g}>0$, which is not possible when the spectrum is symmetric about the origin. In both cases, we use the following two facts from Lemma~\ref{lem:struct}.
\begin{enumerate}
  \item[(a)] If $u\in V(T_i)$, $w\in V(T_j)$ and $\E_{uw}\neq 0$, then $d(v_i,v_j)=k$.
  \item[(b)] If $d(v_i,v_j)=k$, then $\E_{u_iu_j}=2h+k>0$.
\end{enumerate}

\medskip\noindent\textbf{Case 1: $g=2k$ is even.}
The only vertex of $C_g$ at distance $k$ from $v_i$ is $v_{i+k}$. Thus, by (a), if $\E_{uw}\ne0$, then $u\in V(T_i)$ and $w\in V(T_{i+k})$ for some $i$. Let
\[
  X=V(T_0)\cup\cdots\cup V(T_{k-1})\qquad\text{and}\qquad
  Y=V(T_k)\cup\cdots\cup V(T_{2k-1}).
\]
Since $T_i$ and $T_{i+k}$ always lie in different sets, every nonzero entry $\E_{uw}$ has one index in $X$ and the other in $Y$. Let $S$ be the diagonal matrix, indexed by $V(U)$, whose $(u,u)$-entry is $1$ if $u\in X$ and $-1$ if $u\in Y$. Then $S^{-1}=S$ and $(S\E S)_{uw}=-\E_{uw}$ for all $u,w$, that is, $S\E S^{-1}=-\E$. Thus, $\E$ and $-\E$ are similar, and so they have the same eigenvalues with the same multiplicities. Thus, the spectrum of $\E$ is symmetric about the origin.

\medskip\noindent\textbf{Case 2: $g=2k+1$ is odd.}
Assume that the spectrum is symmetric about the origin. Then the eigenvalues split into pairs $\xi,-\xi$ together with some zeros. Since $g$ is odd, $\xi^{g}+(-\xi)^{g}=0$ for every pair, so
\begin{equation}\label{eq:trace}
  \operatorname{tr}\E^{g}=\sum_{i=1}^{n}\xi_i(U)^{g}=0.
\end{equation}
On the other hand, as all entries of $\E$ are nonnegative, the diagonal entry $(\E^{g})_{u_0u_0}$ is the sum of the weights of all closed walks of length $g$ that start and end at $u_0$, where the weight of a walk is the product of the entries of $\E$ along its steps. One of these walks is
\[
  u_0 \to u_k \to u_{2k} \to \cdots \to u_{gk} = u_0 .
\]
Since $d(v_{jk},v_{(j+1)k})=k$ for every $j$, each of its $g$ steps has weight $2h+k$ by (b). Hence the weight of this walk is $(2h+k)^g$. Since all other walks have nonnegative weights, we get
\[
  \operatorname{tr}\E^{g}\ \ge\ (\E^{g})_{u_0u_0}\ \ge\ (2h+k)^{g}\ >\ 0,
\]
which contradicts \eqref{eq:trace}. Hence for odd $g$, the spectrum is not symmetric about the origin.
\end{proof}


\section{Invertibility of the eccentricity matrices of unicyclic graphs}\label{sec:inv}
In this section, we consider the unicyclic graph  $U=C_g(T_0,T_1,\ldots,T_{g-1})$ whose branches $T_i$ may be trivial and establish the conditions when $\E(U)$ is invertible.

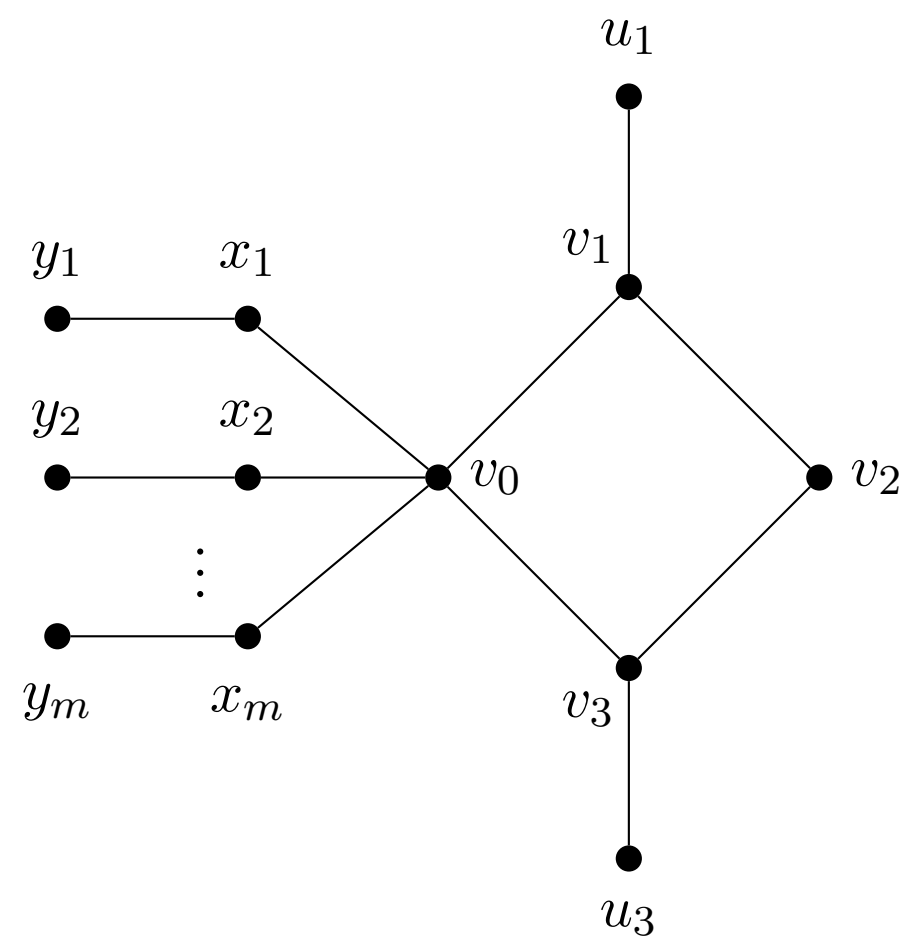
\begin{figure}[htbp]
\centering
\begin{tikzpicture}[every node/.style={circle,draw,inner sep=1.6pt,fill=black},
                    scale=1.05]
  \node[label=right:$v_0$]       (v0) at (0,0)      {};
  \node[label=above left:$v_1$]  (v1) at (1.2,1.2)  {};
  \node[label=right:$v_2$]       (v2) at (2.4,0)    {};
  \node[label=below left:$v_3$]  (v3) at (1.2,-1.2) {};
  \node[label=above:$u_1$]       (a)  at (1.2,2.4)  {};
  \node[label=below:$u_3$]       (c)  at (1.2,-2.4) {};
  \node[label=above:$x_1$] (x1) at (-1.2,1.0)  {};
  \node[label=above:$y_1$] (y1) at (-2.4,1.0)  {};
  \node[label=above:$x_2$] (x2) at (-1.2,0)    {};
  \node[label=above:$y_2$] (y2) at (-2.4,0)    {};
  \node[label=below:$x_m$] (xm) at (-1.2,-1.0) {};
  \node[label=below:$y_m$] (ym) at (-2.4,-1.0) {};
  \draw (v0)--(v1)--(v2)--(v3)--(v0);
  \draw (v1)--(a);  \draw (v3)--(c);
  \draw (v0)--(x1)--(y1);
  \draw (v0)--(x2)--(y2);
  \draw (v0)--(xm)--(ym);
  \node[draw=none,fill=none] at (-1.5,-0.5) {$\vdots$};
\end{tikzpicture}
\caption{The graph $G_m=C_4(S(K_{1,m}),P_2,K_1,P_2)$, which has $2m+6$ vertices.}
\label{fig:Gm}
\end{figure}

First, let us recall some notations, which will be used throughout this section. Let $h_i$ be the height of $T_i$ and $H_i=\max\{d(v_i,w):w\in V(U)\setminus V(T_i)\}$, and $u_i$ is a vertex of $T_i$ with $d(u_i,v_i)=h_i$. Clearly,
\begin{equation}\label{eq:Hi}
  e(v_i)=\max\{h_i,H_i\}\qquad\text{for every } i.
\end{equation}
For $m\ge1$, let $S(K_{1,m})$ be the tree obtained from the star $K_{1,m}$ by subdividing each edge once, rooted at the centre of the star. Let $G_m=C_4(S(K_{1,m}),P_2,K_1,P_2)$, that is, $G_m$ is obtained from the cycle $v_0v_1v_2v_3v_0$ by attaching the $m$ paths $v_0x_ty_t$ $(1\le t\le m)$ at $v_0$, a pendant vertex $u_1$ at $v_1$ and a pendant vertex $u_3$ at $v_3$ (see Figure~\ref{fig:Gm}). The graph $G_m$ has girth $4$ and $2m+6$ vertices.

In the following theorem, we consider the unicyclic graph $U\in\mathcal{U}_{n,g}$ and establish a necessary and sufficient condition for $\E(U)$ to be invertible.


\begin{thm}\label{thm:inv}
Let $U=C_g(T_0,T_1,\ldots,T_{g-1})\in\mathcal{U}_{n,g}$. Then $\E(U)$ is invertible if and only if $T_i\cong P_2$ for every $i$.
\end{thm}

\begin{proof}
By Theorem~\ref{thm:rank-inertia}, it follows that $\rank\E(U)=2g$. Thus, $\E(U)$ is invertible if and only if $n=2g$. Since every $T_i$ is nontrivial, $n=\sum_{i=0}^{g-1}|V(T_i)|\ge2g$, with equality if and only if $|V(T_i)|=2$, that is, $T_i\cong P_2$, for every $i$. This completes the proof.
\end{proof}

Now, we consider a general unicyclic graph $U$ on $n$ vertices and prove a few lemmas, which will be used to prove our main result of this section given in  Theorem~\ref{thm:main-inv}. The following lemma holds for every connected graph.

\begin{lem}\label{lem:twins}
Let $G$ be a connected graph, and let $y$ and $y'$ be distinct vertices of $G$ with $N_G(y)\setminus\{y'\}=N_G(y')\setminus\{y\}$. If $d(y,y')<e(y)$, then $R_y=R_{y'}$, and $\E(G)$ is singular. In particular, this holds if $y$ and $y'$ are pendant vertices with the same neighbour and $e(y)\ge3$.
\end{lem}

\begin{proof}
Let $z\in V(G)\setminus\{y,y'\}$, and let $x$ be the vertex after $y$ on a shortest path from $y$ to $z$. If $x=y'$, then $d(y',z)<d(y,z)$. Otherwise $x$ is a neighbour of $y'$, and $d(y',z)\le1+d(x,z)=d(y,z)$. Thus, $d(y',z)\le d(y,z)$, and in the same way $d(y,z)\le d(y',z)$. Hence $d(y,z)=d(y',z)$ for every $z\in V(G)\setminus\{y,y'\}$. Since $d(y,y')<e(y)$, the eccentricity $e(y)$ is attained at a vertex outside $\{y,y'\}$, and so $e(y')\ge e(y)>d(y,y')$. In the same way $e(y)\ge e(y')$, and hence $e(y)=e(y')$. Thus, $\E(G)_{yz}=\E(G)_{y'z}$ for every $z\in V(G)\setminus\{y,y'\}$. Moreover, $d(y,y')<e(y)=e(y')$ gives $\E(G)_{yy'}=\E(G)_{y'y}=0$. Therefore, $R_y=R_{y'}$. Finally, if $y$ and $y'$ are pendant vertices with the same neighbour, then $N_G(y)=N_G(y')$ and $d(y,y')=2$, and the last assertion follows.
\end{proof}

\begin{lem}\label{lem:h2loc}
Let $U=C_g(T_0,T_1,\ldots,T_{g-1})$ be a unicyclic graph, and let $i\in\{0,1,\ldots,g-1\}$ with $h_i=H_i=2$. Then $d(u,v_i)\le2$ and
\[
  e(u)=d(u,v_i)+2\qquad\text{for every }u\in V(U).
\]
Moreover, for distinct $u,w\in V(U)$,
\[
  \E(U)_{uw}\ne0\iff \max\{d(u,v_i),d(w,v_i)\}=2
  \ \text{ and }\ d(u,w)=d(u,v_i)+d(w,v_i),
\]
and in this case $\E(U)_{uw}=\min\{d(u,v_i),d(w,v_i)\}+2$.
\end{lem}
\begin{proof}
Without loss of generality, assume that $i=0$. By \eqref{eq:Hi}, $e(v_0)=\max\{h_0,H_0\}=2$. Thus, $d(u,v_0)\le2$ and $e(u)\le d(u,v_0)+2$ for every $u$. Let $w_0\notin V(T_0)$ with $d(v_0,w_0)=H_0=2$. If $u\in V(T_0)$, then $d(u,w_0)=d(u,v_0)+2$. If $u\notin V(T_0)$, then $d(u,u_0)=d(u,v_0)+2$ by \eqref{eq:branch-dist}. Hence $e(u)=d(u,v_0)+2$.

Now, let $u$ and $w$ be distinct vertices of $U$. By the above, $\E(U)_{uw}\ne0$ if and only if $d(u,w)=\min\{d(u,v_0),d(w,v_0)\}+2$. Since $d(u,v_0)\le2$ and $d(w,v_0)\le2$,
\begin{align*}
  d(u,w)\ &\le\ d(u,v_0)+d(w,v_0)
  \ =\ \min\{d(u,v_0),d(w,v_0)\}+\max\{d(u,v_0),d(w,v_0)\}\\
  &\le\ \min\{d(u,v_0),d(w,v_0)\}+2 .
\end{align*}
Therefore, $\E(U)_{uw}\ne0$ if and only if both of these inequalities are equalities, that is, if and only if $d(u,w)=d(u,v_0)+d(w,v_0)$ and $\max\{d(u,v_0),d(w,v_0)\}=2$. In this case $\E(U)_{uw}=\min\{d(u,v_0),d(w,v_0)\}+2$.
\end{proof}

\begin{lem}\label{lem:singular}
Let $U=C_g(T_0,T_1,\ldots,T_{g-1})$ be a unicyclic graph, and let $i\in\{0,1,\ldots,g-1\}$. Then the following hold.
\begin{enumerate}
  \item[(a)] If $H_i<h_i$, then $\E(U)$ is singular.
  \item[(b)] Let $H_i\ge h_i$, and let $u$ and $u'$ be distinct vertices of $T_i$ with $d(u,v_i)<h_i$ and $d(u',v_i)<h_i$. If $H_i>h_i$, or if $u$ and $u'$ lie in the same connected component of $T_i-v_i$, then $R_u$ and $R_{u'}$ are multiples of the same vector, and $\E(U)$ is singular.
\end{enumerate}
\end{lem}
\begin{proof}
(a) Let $v,v'\in V(U)\setminus V(T_i)$. Since $d(v_i,v')\le H_i$, the triangle inequality and \eqref{eq:ecc-lower} give
\[
  d(v,v')\le d(v,v_i)+H_i<d(v,v_i)+h_i\le e(v),
\]
and in the same way $d(v,v')<e(v')$. Hence $\E(U)_{vv'}=0$. Let $w$ and $w'$ be the two neighbours of $v_i$ on $C_g$. By \eqref{eq:Hi}, $e(v_i)=h_i$. Thus, $e(w)\le1+e(v_i)=1+h_i$, and \eqref{eq:ecc-lower} gives $e(w)=1+h_i$. In the same way, $e(w')=1+h_i$. Also, $d(w,x)=1+d(v_i,x)=d(w',x)$ for every $x\in V(T_i)$. Thus, $\E(U)_{wx}=\E(U)_{w'x}$ for $x\in V(T_i)$, and $\E(U)_{wx}=\E(U)_{w'x}=0$ for $x\notin V(T_i)$. Therefore, $R_w=R_{w'}$, and $\E(U)$ is singular.

(b) Since $H_i\ge h_i$, \eqref{eq:Hi} gives $e(v_i)=H_i=d(v_i,v)$ for some $v\notin V(T_i)$. Thus, by Lemma~\ref{lem:tree}, $e(x)=d(x,v_i)+H_i$ for every $x\in V(T_i)$. Now, let $x\in\{u,u'\}$.

Let $w\notin V(T_i)$. Since every path from $x$ to $w$ passes through $v_i$, \eqref{eq:ecc-lower} gives $e(w)\ge d(w,v_i)+h_i>d(w,v_i)+d(v_i,x)=d(x,w)$. Thus, $\E(U)_{xw}\ne0$ if and only if $d(x,w)=e(x)$, that is, $d(v_i,w)=H_i$.

Let $w\in V(T_i)$. Since $d(x,v_i)<h_i$ and $d(w,v_i)\le h_i$,
\begin{align*}
  d(x,w)&\le d(x,v_i)+d(w,v_i)\le\min\{d(x,v_i),d(w,v_i)\}+h_i\\
  &\le\min\{d(x,v_i),d(w,v_i)\}+H_i=\min\{e(x),e(w)\}.
\end{align*}
Thus, $\E(U)_{xw}\ne0$ if and only if $H_i=h_i$, $d(w,v_i)=h_i$ and $d(x,w)=d(x,v_i)+d(v_i,w)$. If $H_i>h_i$, then no such $w$ exists. If $H_i=h_i$, then $u$ and $u'$ lie in the same connected component $B$ of $T_i-v_i$, and the last condition means that $w\notin V(B)$.

In all cases, the condition for $\E(U)_{xw}\ne0$ does not depend on the choice of $x\in\{u,u'\}$, and then $\E(U)_{xw}=d(x,w)=e(x)$. Hence $R_u=e(u)\chi_X^T$ and $R_{u'}=e(u')\chi_X^T$ for the same set $X\subseteq V(U)$. Therefore, $\E(U)$ is singular.
\end{proof}

\begin{lem}\label{prop:common}
Let $U=C_g(T_0,T_1,\ldots,T_{g-1})$ be a unicyclic graph such that $\E(U)$ is invertible. Then the following hold.
\begin{enumerate}
  \item[(i)] $h_i\le H_i$ and $h_i\le2$ for every $i$.
  \item[(ii)] If $h_i=2$ for some $i$, then $H_i=2$ and $g=4$.
  \item[(iii)] If $h_i\le1$ for some $i$, then $T_i$ is trivial or $T_i\cong P_2$.
\end{enumerate}
\end{lem}

\begin{proof}
Parts (i) and (ii) follow from Lemma~\ref{lem:singular}, together with Lemma~\ref{lem:h2loc} when $g$ is odd. Part (iii) follows from Lemma~\ref{lem:twins}.

(i) By Lemma~\ref{lem:singular}(a), $h_i\le H_i$ for every $i$. Assume, on the contrary, that $h_i\ge3$ for some $i$, and let $v_ix_1x_2\cdots x_{h_i}$ be a path in $T_i$. Then $x_1$ and $x_2$ lie in the same connected component of $T_i-v_i$, and $d(x_1,v_i)=1<d(x_2,v_i)=2<h_i$. Thus, $\E(U)$ is singular by Lemma~\ref{lem:singular}(b), a contradiction. Hence $h_i\le2$ for every $i$.

(ii) Without loss of generality, assume that $i=0$. By (i), $H_0\ge2$. Assume, on the contrary, that $H_0>2$, and let $x$ be a child of $v_0$. Since $d(v_0,v_0)=0<2$ and $d(x,v_0)=1<2$, Lemma~\ref{lem:singular}(b) shows that $\E(U)$ is singular, a contradiction. Hence $H_0=2$, and Lemma~\ref{lem:h2loc} applies with $i=0$. Since $v_j\notin V(T_0)$ for every $j\ne0$, we have $k\le H_0=2$, and so $g\le5$.

Assume, on the contrary, that $g$ is odd. Then $g=2k+1$. Let $p=v_k$ and $q=v_{k+1}$. Then $d(p,v_0)=d(q,v_0)=k$ and $d(p,q)=1$. We show that $R_p=R_q$. Let $x\in\{p,q\}$, and let $w\in V(U)\setminus\{p,q\}$. If $w\in V(T_0)$, then $d(x,w)=d(x,v_0)+d(v_0,w)$ by \eqref{eq:branch-dist}. If $w\in V(T_j)$ for some $j\ne0$, then
\[
  d(x,w)\le d(x,v_j)+d(v_j,w)\le k+d(v_j,w)<k+d(v_0,v_j)+d(v_j,w)=d(x,v_0)+d(v_0,w).
\]
Thus, by Lemma~\ref{lem:h2loc}, $\E(U)_{xw}\ne0$ if and only if $w\in V(T_0)$ and $\max\{k,d(w,v_0)\}=2$, and in this case $\E(U)_{xw}=\min\{k,d(w,v_0)\}+2$. This does not depend on the choice of $x\in\{p,q\}$. Moreover, $d(p,q)=1<k+2=\min\{e(p),e(q)\}$, and so $\E(U)_{pq}=\E(U)_{qp}=0$. Hence $R_p=R_q$, and $\E(U)$ is singular, a contradiction. Therefore, $g$ is even, and $g=4$.

(iii) Since $h_i\le1$, we have $T_i\cong K_{1,s}$, where $s=\deg_{T_i}(v_i)$, and it is enough to show that $s\le1$. Assume, on the contrary, that $s\ge2$, and let $y$ and $y'$ be two leaves of $T_i$. They are pendant vertices of $U$ with the common neighbour $v_i$. Since $d(y,w)=1+d(v_i,w)$ for every $w\notin V(T_i)$, we have $e(y)\ge1+H_i$. If $g\ge4$, then $H_i\ge k\ge2$. If $g=3$ and $h_j\ge1$ for some $j\ne i$, then $H_i\ge d(v_i,v_j)+h_j\ge2$. In both cases, $e(y)\ge3$, and $\E(U)$ is singular by Lemma~\ref{lem:twins}, a contradiction. In the remaining case, $g=3$, and $T_{i+1}$ and $T_{i+2}$ are trivial. Then $N_U(v_{i+1})\setminus\{v_{i+2}\}=\{v_i\}=N_U(v_{i+2})\setminus\{v_{i+1}\}$ and $e(v_{i+1})\ge d(v_{i+1},y)=2>d(v_{i+1},v_{i+2})$. Thus, $\E(U)$ is singular by Lemma~\ref{lem:twins}, again a contradiction. Hence $s\le1$, that is, $T_i$ is trivial or $T_i\cong P_2$.
\end{proof}

Using the previous lemmas, in the following theorem, we establish a necessary and sufficient condition for the invertibility of eccentricity matrix of a unicyclic graph $U=C_g(T_0,T_1,\ldots,T_{g-1})$ with $h_i\ge2$ for some $i$. Moreover, we obtain a necessary condition for the invertibility of $\mathcal{E}(U)$ when $h_i\le1$ for every $i$.

\begin{thm}\label{thm:main-inv}
Let $U=C_g(T_0,T_1,\ldots,T_{g-1})$ be a unicyclic graph on $n$ vertices.
\begin{enumerate}
  \item[(a)] If $h_i\ge2$ for some $i$, then $\E(U)$ is invertible if and only if $U\cong G_m$ for some $m\ge1$. In this case, $g=4$ and $n=2m+6$.
  \item[(b)] If $h_i\le1$ for every $i$ and $\E(U)$ is invertible, then every branch of $U$ is trivial or isomorphic to $P_2$.
\end{enumerate}
\end{thm}

\begin{proof}
We first prove two claims about the graphs $G_m$, and we use them to prove (a). Then we prove (b).

\smallskip
\noindent\textbf{Claim 1. If $\E(U)$ is invertible and $h_i=2$ for some $i$, then $U\cong G_m$, where $m=\deg_{T_i}(v_i)$.}

\noindent\emph{Proof of Claim~1.} By Lemma~\ref{prop:common}(ii), $H_i=2$ and $g=4$. Without loss of generality, assume that $i=0$. Let $X=\{u\in V(U):d(u,v_0)\le1\}$ and $Y=\{u\in V(U):d(u,v_0)=2\}$. By Lemma~\ref{lem:h2loc}, $d(u,v_0)\le2$ for every $u$, and so $V(U)=X\cup Y$. Also by Lemma~\ref{lem:h2loc}, $\E(U)_{v_0w}=2$ if $w\in Y$, and $\E(U)_{v_0w}=0$ otherwise. Thus, $R_{v_0}=2\chi_Y^{T}$.

First, we show that $T_0\cong S(K_{1,m})$. Assume, on the contrary, that some child $x$ of $v_0$ has no child. Then $v_0$ is the only neighbour of $x$. Hence $d(x,w)=1+d(v_0,w)=3=d(x,v_0)+d(w,v_0)$ for every $w\in Y$, and Lemma~\ref{lem:h2loc} gives $\E(U)_{xw}=\min\{1,2\}+2=3$ for $w\in Y$, and $\E(U)_{xw}=0$ otherwise. Thus, $R_x=3\chi_Y^{T}=\tfrac32R_{v_0}$, and $\E(U)$ is singular, a contradiction. Hence every child $x$ of $v_0$ has a child. Since $h_0=2$, the children of $x$ are pendant vertices of $U$. If $x$ had two children $y$ and $y'$, then $y$ and $y'$ would be pendant vertices with the common neighbour $x$, and $e(y)=d(y,v_0)+2=4$ by Lemma~\ref{lem:h2loc}. Then $\E(U)$ would be singular by Lemma~\ref{lem:twins}. Therefore, every child of $v_0$ has exactly one child, that is, $T_0\cong S(K_{1,m})$ with $m=\deg_{T_0}(v_0)\ge1$.

Next, we find the other branches. Since $d(v_2,v_0)=2$, the tree $T_2$ is trivial, and since $d(v_1,v_0)=d(v_3,v_0)=1$, we get $h_1,h_3\le1$. Thus, $T_1$ and $T_3$ are trivial or isomorphic to $P_2$ by Lemma~\ref{prop:common}(iii). For $j\in\{1,3\}$, let $s_j=|V(T_j)|-1\in\{0,1\}$.

By Lemma~\ref{lem:h2loc}, $\E(U)_{uw}=0$ for all $u,w\in X$. Thus, the $|X|$ rows $R_u$, $u\in X$, have all their nonzero entries in the $|Y|$ columns indexed by $Y$. Since $\E(U)$ is invertible, these rows are linearly independent, and hence $|X|\le|Y|$. Let $v_0x_ty_t$ $(1\le t\le m)$ be the paths in $T_0$ from $v_0$ to the leaves of $T_0$. Then $X=\{v_0,v_1,v_3,x_1,\ldots,x_m\}$, and $Y$ consists of $v_2,y_1,\ldots,y_m$ and the $s_1+s_3$ leaves of $T_1$ and $T_3$. Thus, $m+3=|X|\le|Y|=m+1+s_1+s_3$. This gives $s_1+s_3\ge2$, and hence $s_1=s_3=1$, that is, $T_1\cong T_3\cong P_2$. Therefore, $U\cong G_m$.

\smallskip
\noindent\textbf{Claim 2. For every $m\ge1$, the matrix $\E(G_m)$ is invertible.}

\noindent\emph{Proof of Claim~2.} Let $U=G_m$. Here $h_0=H_0=2$, so Lemma~\ref{lem:h2loc} applies with $i=0$. By Figure~\ref{fig:Gm}, the vertices split into the sets
\begin{align*}
  X&=\{u:d(u,v_0)\le1\}=\{v_0,v_1,v_3,x_1,\ldots,x_m\},\\
  Y&=\{u:d(u,v_0)=2\}=\{v_2,u_1,u_3,y_1,\ldots,y_m\},
\end{align*}
and both sets have $m+3$ vertices. By Lemma~\ref{lem:h2loc}, $\E(U)_{uw}=0$ for all $u,w\in X$. Thus, with the vertices in the order $X$ followed by $Y$,
\[
  \E(U)=\begin{pmatrix}0&B\\ B^{T}&L\end{pmatrix}
\]
for some matrices $B$ and $L$ of order $m+3$. Moving the last $m+3$ columns of $\E(U)$ in front of the first $m+3$ columns needs $(m+3)^{2}$ interchanges of columns, and this move gives a block lower triangular matrix. Thus, $\det\E(U)=(-1)^{(m+3)^{2}}\det B\det B^{T}=(-1)^{m+1}(\det B)^{2}$, and it is enough to show that $B$ is nonsingular.

By Lemma~\ref{lem:h2loc}, for $u\in X$ and $w\in Y$, the entry $\E(U)_{uw}$ is $d(u,v_0)+2$ if $d(u,w)=d(u,v_0)+2$, and it is $0$ otherwise. Since $d(v_0,w)=2$ for every $w\in Y$, all entries in the row of $v_0$ in $B$ are equal to $2$. For the other rows, we have $d(u,v_0)=1$, so the entry is $3$ exactly when $d(u,w)=3$. By Figure~\ref{fig:Gm}, this happens for
\begin{align*}
  w&\in\{u_3,y_1,\ldots,y_m\} &&\text{if } u=v_1,\\
  w&\in\{u_1,y_1,\ldots,y_m\} &&\text{if } u=v_3,\\
  w&\in\{v_2,u_1,u_3\}\cup\{y_t:t\ne s\} &&\text{if } u=x_s.
\end{align*}
Now, let $\beta=(\beta_w)_{w\in Y}$ with $B\beta=0$, and let $\sigma=\sum_{w\in Y}\beta_w$. The row of $v_0$ gives $2\sigma=0$, so $\sigma=0$. For $1\le s\le m$, the row of $x_s$ gives $3(\sigma-\beta_{y_s})=0$, so $\beta_{y_s}=0$. Then the rows of $v_1$ and $v_3$ give $\beta_{u_3}=0$ and $\beta_{u_1}=0$, respectively, and $\sigma=0$ gives $\beta_{v_2}=0$. Thus, $\beta=0$, and $B$ is nonsingular.

\smallskip
\noindent\emph{Proof of (a).} Let $h_i\ge2$ for some $i$. If $\E(U)$ is invertible, then $h_i=2$ by Lemma~\ref{prop:common}(i), and $U\cong G_m$ for some $m\ge1$ by Claim~1. Conversely, if $U\cong G_m$, then $\E(U)$ is invertible by Claim~2. Since $G_m$ has girth $4$ and $2m+6$ vertices, this proves (a).

\smallskip
\noindent\emph{Proof of (b).} Let $h_i\le1$ for every $i$, and let $\E(U)$ be invertible. Then every branch of $U$ is trivial or isomorphic to $P_2$ by Lemma~\ref{prop:common}(iii).
\end{proof}

\begin{rmk}\label{rmk:not-sufficient}
The converse of Theorem~\ref{thm:main-inv}(b) is not true. For example, every branch of the graph $U=C_3(P_2,P_2,K_1)$ in Figure~\ref{fig:C3} is trivial or isomorphic to $P_2$, but $\E(U)$ is not invertible, since $R_{v_2}=R_{v_0}+R_{v_1}$.

\begin{figure}[htbp]
\centering
\begin{tikzpicture}[every node/.style={circle,draw,inner sep=1.6pt,fill=black},
                    scale=1.05]
  \node[label=below:$v_0$] (v0) at (0,0)     {};
  \node[label=below:$v_1$] (v1) at (1.6,0)   {};
  \node[label=above:$v_2$] (v2) at (0.8,1.3) {};
  \node[label=below:$u_0$] (u0) at (-1.4,0)  {};
  \node[label=below:$u_1$] (u1) at (3.0,0)   {};
  \draw (v0)--(v1)--(v2)--(v0);
  \draw (u0)--(v0);  \draw (v1)--(u1);
\end{tikzpicture}
\caption{The graph $C_3(P_2,P_2,K_1)$.}
\label{fig:C3}
\end{figure}
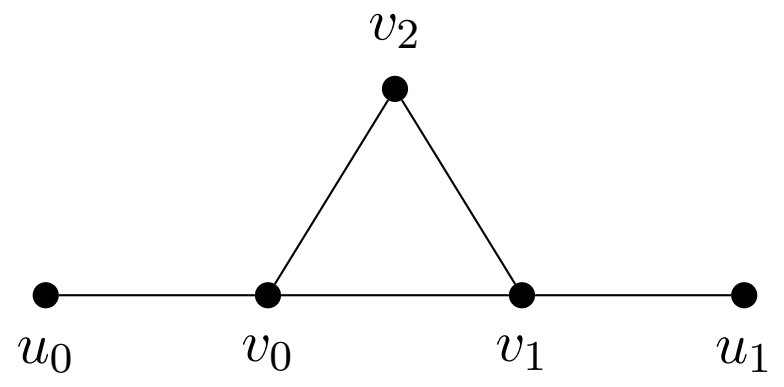
\end{rmk}



\FloatBarrier

\bibliographystyle{abbrv}
\bibliography{ref}
\end{document}